\documentclass{gtpart}     
\usepackage[dvipdfmx]{graphicx}  
\usepackage{pinlabel}  
\usepackage[all]{xy}
\usepackage{amsmath,amssymb,amscd}
\title{Coverage Criterion in Sensor Networks Stable under \\ Perturbation}

\author{Yasuaki Hiraoka} 
\givenname{Yasuaki}
\surname{Hiraoka}
\address{Institute of Mathematics for Industry\\Kyushu University\\
              744, Motooka, Nishi-ku, Fukuoka, 819-0395, Japan}
\email{hiraoka@imi.kyushu-u.ac.jp}
\urladdr{}

\author{Genki Kusano} 
\givenname{Genki}
\surname{Kusano}
\address{Graduate School of Mathematics\\Kyushu University\\
              744, Motooka, Nishi-ku, Fukuoka, 819-0395, Japan}
\email{ma214017@math.kyushu-u.ac.jp}

\arxivreference{1409.7483}
\arxivpassword{4ktgm}

\volumenumber{}
\issuenumber{}
\publicationyear{}
\papernumber{}
\startpage{}
\endpage{}
\doi{}
\MR{}
\Zbl{}
\received{}
\revised{}
\accepted{}
\published{}
\publishedonline{}
\proposed{}
\seconded{}
\corresponding{}
\editor{}
\version{}

\newtheorem{thm}{Theorem}[section]   
\newtheorem{lem}[thm]{Lemma}         
\newtheorem{prop}[thm]{Proposition}
\newtheorem{cor}[thm]{Corollary}
\newtheorem{exa}[thm]{Example}
\newtheorem*{dSG}{dSG Criterion}
\newtheorem*{main}{Main theorem}

\makeatletter 
\@tempcnta\z@
\loop\ifnum\@tempcnta<26
\advance\@tempcnta\@ne
\expandafter\edef\csname f\@Alph\@tempcnta\endcsname{\noexpand\mathfrak{\@Alph\@tempcnta}}
\expandafter\edef\csname g\@alph\@tempcnta\endcsname{\noexpand\mathfrak{\@alph\@tempcnta}}
\expandafter\edef\csname l\@Alph\@tempcnta\endcsname{\noexpand\mathbb{\@Alph\@tempcnta}}
\expandafter\edef\csname c\@Alph\@tempcnta\endcsname{\noexpand\mathcal{\@Alph\@tempcnta}}
\repeat

\newcommand{\ee}{{\varepsilon}}
\newcommand{\pa}{{\partial}}
\newcommand{\one}{\mbox{1}\hspace{-0.25em}\mbox{l}}

\newcommand{\ra}{{\rightarrow}}

\newcommand{\rra}{{\rightrightarrows}}
\newcommand{\hra}{{\hookrightarrow}}

\newcommand{\dis}{{\rm dis}}
\newcommand{\Hom}{{\rm Hom}}
\newcommand{\Nr}{{N_{\hat{r}}(\pa D)}}
\newcommand{\im}{\mathop{\rm im}\nolimits}

\providecommand{\abs}[1]{\left\lvert#1\right\rvert}
\providecommand{\norm}[1]{\left\lVert#1\right\rVert}

\begin{document}

\begin{abstract}    
To the coverage problem of sensor networks, Vin de Silva and Robert Ghrist developed several approaches \cite{dSG1,dSG2} based on (persistent) homology theory. 
Their criteria for the coverage are formulated on the Rips complexes constructed by the sensors, in which their locations are supposed to be fixed. 
However, the sensors are in general affected by perturbations (e.g., natural phenomena), and hence the stability of the coverage criteria should be also discussed. 
In this paper, we present a coverage theorem stable under perturbation. 
Furthermore, we also introduce a method of eliminating redundant cover after perturbation.
The coverage theorem is derived by extending the Rips interleaving theorem studied in the paper \cite{CSO} into an appropriate relative version. 

\end{abstract}

\maketitle


\section{INTRODUCTION}

\begin{figure}[htbp]
 \begin{minipage}{0.48\hsize}
 \begin{center}
  \includegraphics[width=60mm]{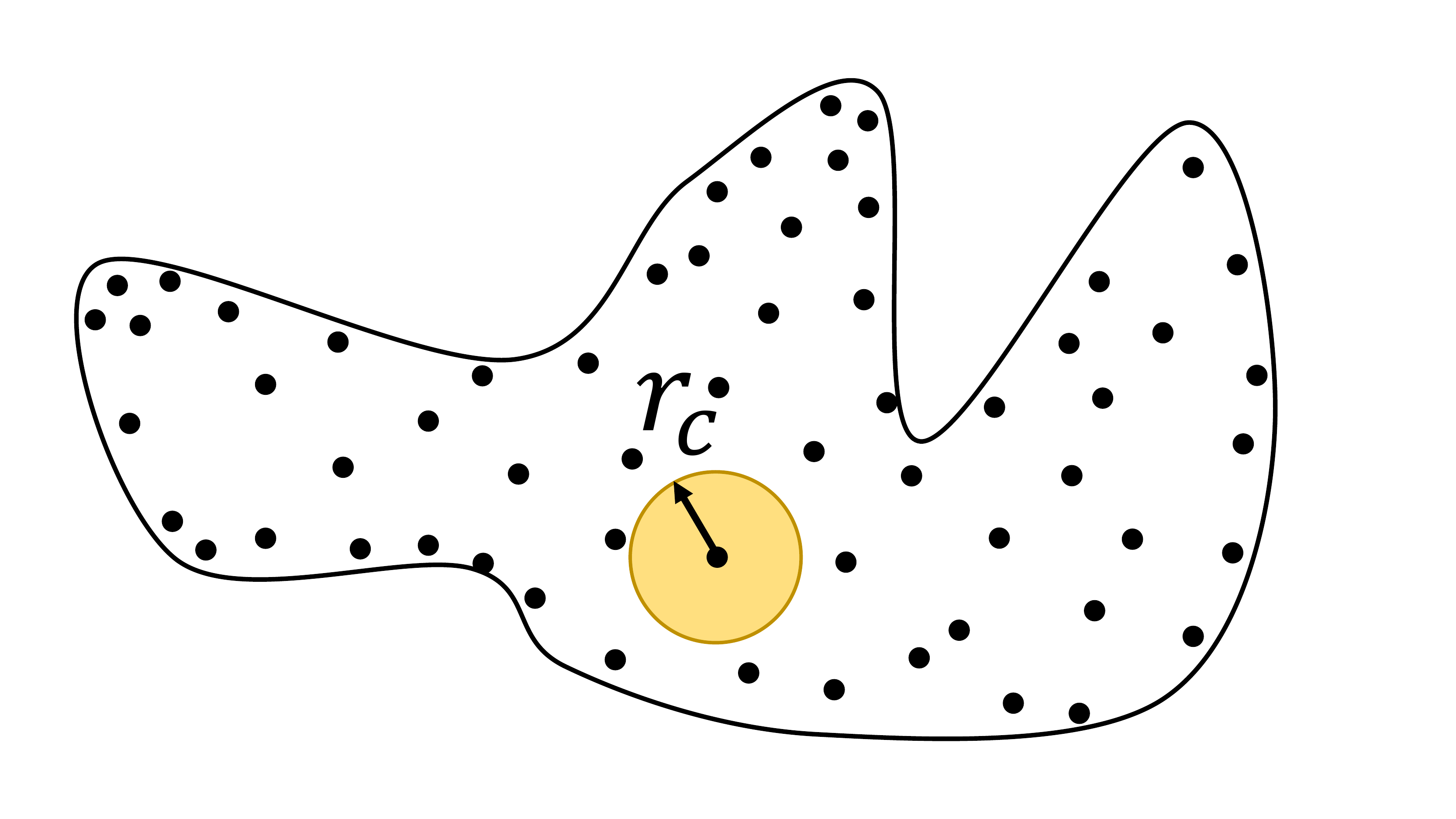}
 \end{center}
\caption{The region surrounded by the curve is a target domain $D$ and the points represent the sensors in $D$. The yellow disk with the radius $r_c$ is the sensing region.}
\label{fig:radi}
 \end{minipage}
 \hspace{4mm}
 \begin{minipage}{0.48\hsize}
 \begin{center}
  \includegraphics[width=60mm]{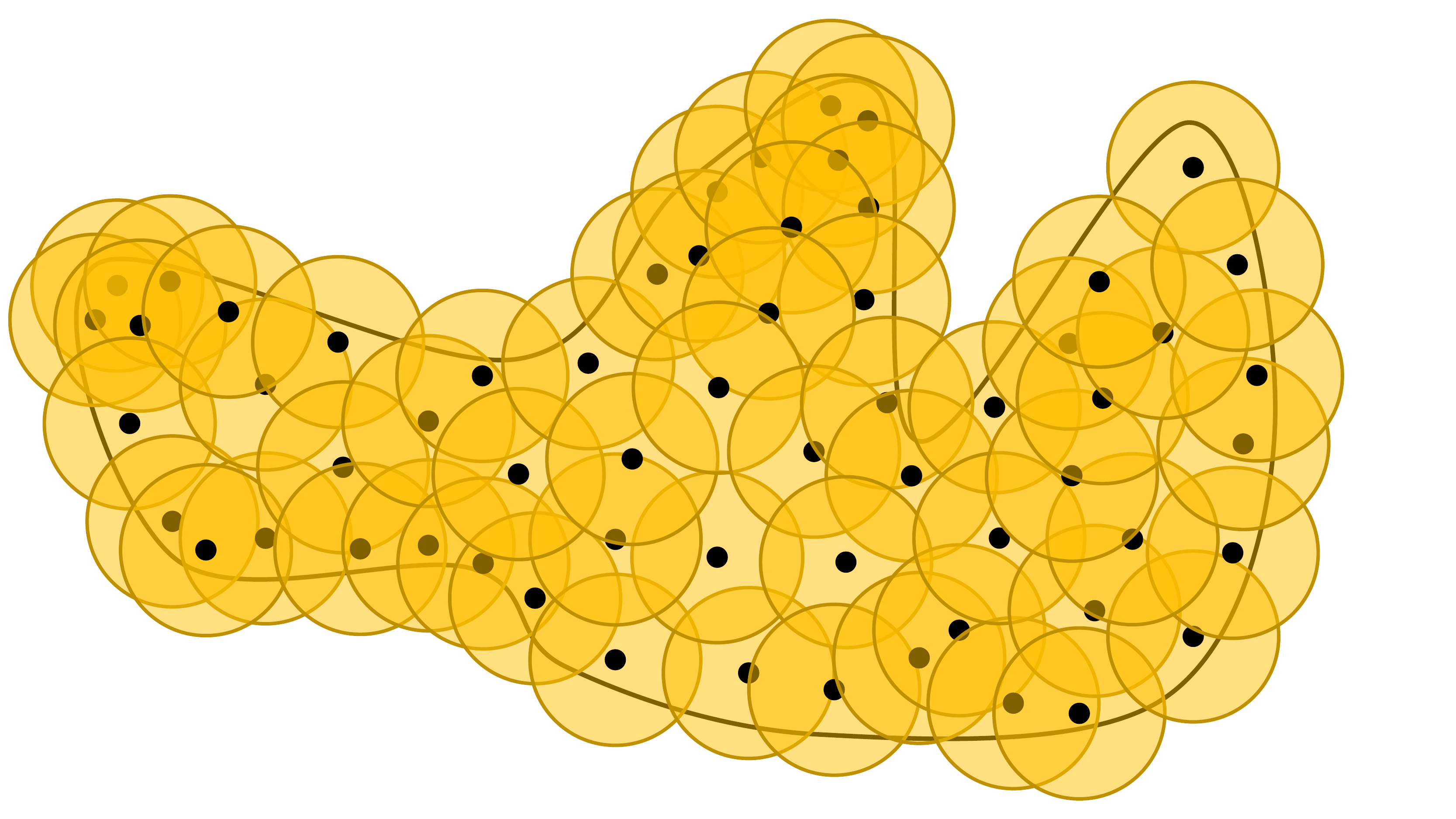}
 \end{center}
\caption{The whole sensing region $\cU(X)$ covers the target domain $D$.}
\label{fig:cover}
 \end{minipage}
\end{figure}

A fundamental task of sensor networks is to extract information about a target domain by using sensors.
Each sensor gathers data around its location, as shown in Figure \ref{fig:radi}.
Here, we are interested in a coverage problem that appears in a variety of settings, such as military, environmental sensing, security, and so on \cite{MA,MK}.
The problem is to determine whether the target domain is covered by the sensing region.
Figure \ref{fig:cover} shows a sensing region covering its target domain.

We define the target domain as a subset $D$ of $\R^d$ and the set of sensors as a finite subset $X=\left\{x_i \in D \right\}$.
Each sensor $x_i$ can monitor its surroundings within a cover radius $r_c$.
Let $B_{r_c}(x_i)=\left\{y \in \R^d \ \middle| \ \|x_i-y\| \leq r_c\right\}$, where $\| \cdot \|$ represents the Euclidean norm on $\R^d$, be the sensing region of $x_i$ and $\cU(X)=\bigcup_{i} B_{r_c}(x_i)$ be the whole sensing region of $X$.
Then the coverage problem is formulated as ``Is the target domain $D$ covered by the whole sensing region $\cU(X)$?''.

From now on, we assume that each sensor does not have ability to gather absolute positional information, e.g., GPS.
In other words, we cannot obtain the coordinates of the sensors.
However, we assume that each sensor can communicate with other sensors if the distance between them is less than a certain communication radius.
Here, we consider the coverage of a restricted target domain $D' \subset D$, where $D'$ is relatively large in $D$.
In this paper, our assumptions are the same as those used in \cite{dSG2}.

\begin{description}

\item[A1.] The cover radius is $r_c$.
\item[A2.] We have two communication radii $r_w$ and $r_s$.
\item[A3.] The communication radii $r_w$, $r_s$ and the cover radius $r_c$ satisfy
\[
r_c \geq \frac{r_s}{\sqrt{2}}, \quad r_w \geq r_s \sqrt{10}. 
\]
\item[A4.] The target domain $D$ is a compact subset of $\R^d$. In addition, the set 
\[
N_f(\pa D):=\left\{x \in D \ \middle| \ \|x-\pa D\| \leq r_f \right\}
\]
defines the set $F$ of fence sensors as $X \cap N_f(\pa D)$, where $r_f$ is the fence detection radius.
\item[A5.] The restricted domain $D-\Nr$ is connected, where
\[
\Nr:=\left\{x \in D \ \middle| \ \|x-\pa D\| \leq \hat{r}:=r_f+\frac{r_s}{\sqrt{2}} \right\} .
\]
\item[A6.] The fence detection hypersurface $\Sigma=\left\{ x \in D \ \middle| \ \|x-\pa D\| = r_f \right\}$ has internal injectivity radius at least $\frac{r_s}{\sqrt{2}}$ and external injectivity radius at least $r_s$.
\end{description}

In these settings, Vin de Silva and Robert Ghrist proposed a criterion to solve the coverage problem from the communication data of sensors. 
Throughout, we call the criterion in \cite[Theorem 3.4]{dSG2} as the dSG criterion for short.

Let $\cR(X;a)$ denote the Rips complex of $X$ with parameter $a$, which we define in Section \ref{preliminary}.
When $F$ is a subset of $X$, $\cR(F;a)$ is a subcomplex of $\cR(X;a)$, and $H_k(\cR(X,F;a))$ means the relative homology group of this pair of the Rips complexes.
\begin{dSG}[Theorem 3.4 of \cite{dSG2}]
\label{dSG}
Let $X$ be a set of sensors in a target domain $D \subset \R^d$ satisfying assumptions \textbf{A1}-\textbf{A6}. 
If the homomorphism
\[
\iota_*:H_d(\cR(X,F;r_s)) \ra H_d(\cR(X,F;r_w))
\]
induced by the inclusion $\iota:\cR(X;r_s) \hra \cR(X;r_w)$ is nonzero,
then the cover $\cU(X)$ of the sensors contains the restricted domain $D-\Nr$.
\end{dSG}

\begin{figure}[htbp]
 \begin{minipage}{0.48\hsize}
 \begin{center}
  \includegraphics[width=60mm]{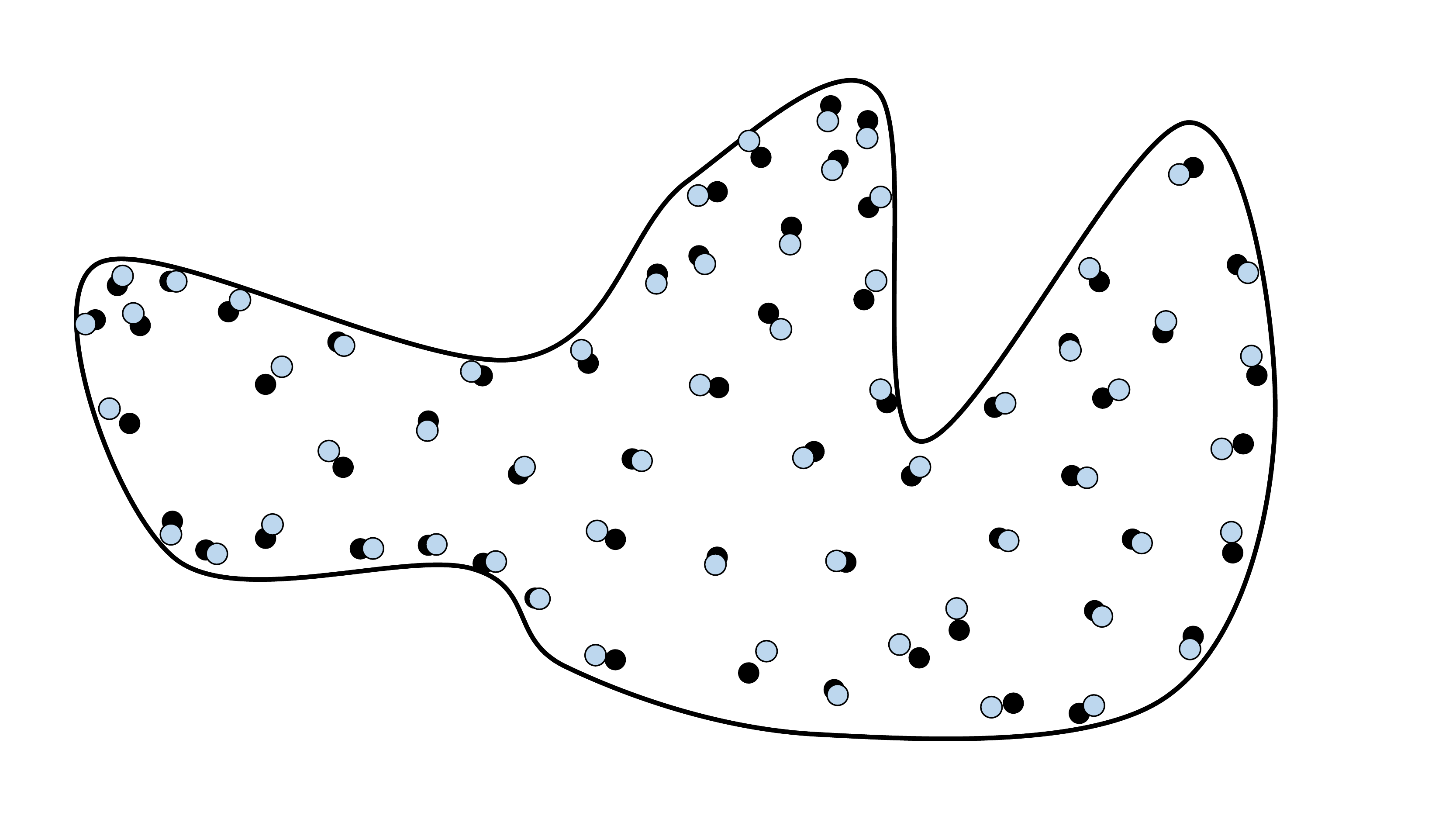}
 \end{center}
\caption{Sensors move to the blue points as a result of a perturbation.}
\label{fig:perturb}
 \end{minipage}
 \hspace{4mm}
 \begin{minipage}{0.48\hsize}
 \begin{center}
  \includegraphics[width=60mm]{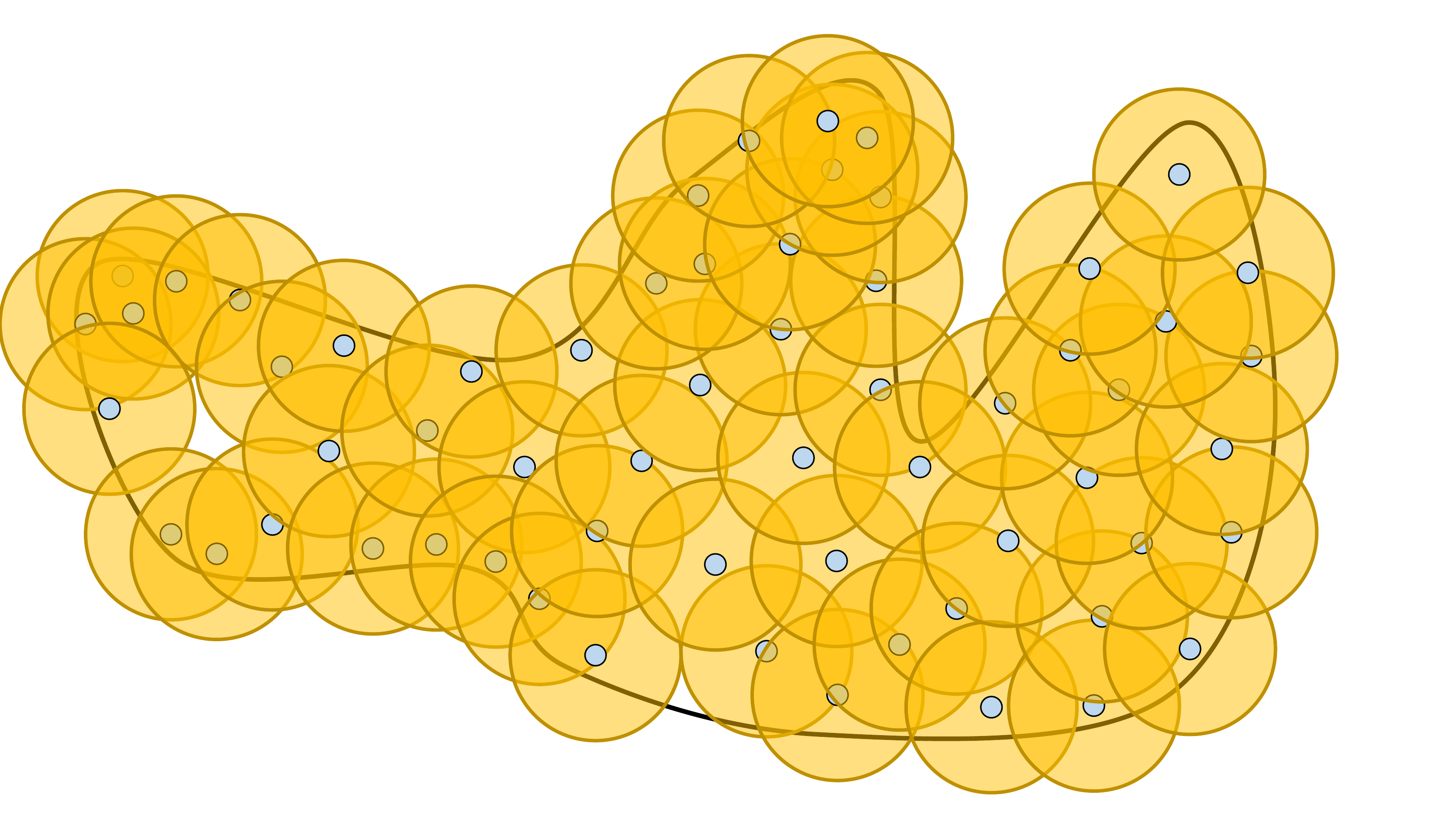}
 \end{center}
\caption{The sensors affected by the perturbation do not cover the target domain, although those located at the original position cover it.}
\label{fig:pcover}
 \end{minipage}
\end{figure}

The dSG criterion applies only to fixed sensors.
However, sensors are not stationary in general.
For example, sensors can be affected by perturbations such as wind or earthquakes.
They may move (Figure \ref{fig:perturb}) and thus the coverage may change (Figure \ref{fig:pcover}).

One option to this situation is to check the dSG criterion again to determine the coverage of the sensors after the perturbation.
Instead, we improve the dSG criterion to be stable under perturbation.
Namely, once we examine the original induced map, automatically we also know the coverage of the sensors even after the perturbation.

\begin{main}
\label{mainthm}
Let $X$ be a set of sensors which satisfies assumptions \textbf{A1}-\textbf{A6} and let $f:X \ra D$ satisfy
$\sup_{x \in X} \|x-f(x)\| \leq \frac{\ee}{2}$, $\ f(F) \subset N_f(\pa D)$, and $f(X-F) \subset D-N_f (\pa D)$.
If the homomorphism
\[
\iota^{\ee}_*:H_d(\cR(X,F;r_s-\ee)) \ra H_d(\cR(X,F;r_w+\ee))
\]
induced by the inclusion $\iota^{\ee}:\cR(X;r_s-\ee) \hra \cR(X;r_w+\ee)$ is nonzero,
then the cover $\cU(f(X))$ of the sensors after the perturbation contains the restricted domain $D-\Nr$.
\end{main}

Here, we view $f$ as a perturbation;
that is, each sensor $x_i$ moves up to a distance $\frac{\ee}{2}$ from its original position.
We emphasize that this criterion does not require the new communication data of the sensors after the perturbation.
In other words, for $X$ and $f$ satisfying the conditions and hypothesis of the main theorem,
we can conclude that the map
\[
j_*:H_d(\cR(f(X),f(F);r_s)) \ra H_d(\cR(f(X),f(F);r_w))
\]
is nonzero, which guarantees the dSG criterion for the perturbed sensors $f(X)$.

Furthermore, if $\alpha \in H_d(\cR(X,F;r_s-\ee)) $ satisfies $\iota^{\ee}_* \alpha \neq 0$, the proof of the main theorem ensures that $f(\alpha) \in H_d(\cR(f(X),f(F);r_s))$ also satisfies $j_* (f(\alpha)) \neq 0$.
This generator $f(\alpha)$ enables us to study how to eliminate redundant cover.
By turning off some sensors so that the cover of the remaining sensors still contains the restricted domain, a smaller cover of $\cU(f(X))$ can contribute to power conservation of active sensors \cite{dSG1}.
Moreover, by applying methods in optimal homologous cycles \cite{DHK, EHi}, an even smaller number of sensors can contain the restricted domain.

In a similar spirit to \cite{dSG1,dSG2}, we here give an example satisfying the assumptions $f(F) \subset N_f(\pa D)$ and $f(X-F) \subset D-N_f(\pa D)$, in which the sensors close to the boundary are controlled.
Let the neighboring sensors
\[
F'=\{x \in X \mid \|x - \pa D\| \leq r_f + \frac{\ee}{2} \}
\]
of $F$ be fixed.
For any sensor $x$ located in $X-F'$, the corresponding perturbed sensor $f(x)$ is in $D-N_f(\pa D)$ because $x$ moves up to a distance $\frac{\ee}{2}$.
Since the other sensors are unchanged, $f(F) \subset N_f(\pa D)$ and $f(X-F) \subset D-N_f(\pa D)$ are satisfied.

A key mathematical concept for studying the stability of the coverage is the {\em interleaving} of persistence modules \cite{CCGGY,CSO}. 
This characterizes the similarity of two persistence modules, and hence some features in one persistence module can be studied by another. 
In this paper, we study the stability property of the coverages between the original and perturbed sensors by interleaving. 
In the derivation of the main theorem, our mathematical contribution is to extend the Rips interleaving theorem studied in \cite{CSO} into an appropriate relative version. 
We also emphasize that the concept of interleavings allows us to compute optimal homologous cycles in the perturbed complex without needing to recompute homology.
That is, we obtain a generator in the unperturbed complex, map it by $f$, then optimize it.

This paper is organized as follows. 
In Section \ref{preliminary}, the basic concepts of correspondences, the Gromov-Hausdorff distance, persistence modules, Rips complexes, Rips interleaving and optimal homologous cycles are introduced.
In Section \ref{relative}, we extend some properties of persistence modules to the relative version and prove relative Rips interleaving.
In Section \ref{power}, some of recent results about eliminating redundant cover are reviewed for the dSG criterion.
In Section \ref{perturbed sensors},we prove the main theorem using relative Rips interleaving and discuss eliminating redundant cover of perturbed sensors.

\section{PRELIMINARIES}
\label{preliminary}
In this section, we review the basic concepts of correspondences, persistence modules and optimal homologous cycles.
We refer the reader to \cite{BBI} for correspondences, \cite{Ha,Mu} for homology groups, \cite{CCGGY,EH} for persistence modules and \cite{DHK, EHi} for optimal homologous cycles.
We basically follow the exposition in \cite{CSO,DHK} with appropriate modifications.
We assume that all vector spaces are defined over arbitrary field $\lF$ and that the coefficient group of homology groups is $\lF$.
In considering optimal homologous cycles, we let the coefficient group of homology groups be $\lR$.

\subsection{Correspondences}
Let $X$ and $Y$ be two sets.
A \textbf{correspondence} $C$ from $X$ to $Y$, denoted by $C:X \rra Y$, is a subset of $X \times Y$ satisfying the following conditions:
for every $x \in X$ there exists at least one $y \in Y$ such that $(x,y) \in C$, and for every $y \in Y$ there exists at least one $x \in X$ such that $(x,y) \in C$.
Since $C=X \times Y$ is a correspondence, there exists at least one correspondence between any two sets.

\begin{exa}
\label{graph}
For any map $f:X \ra Y$, its graph $G(f)=\left\{(x,f(x)) \in X \times Y \ \middle| \ x \in X\right\} \subset X \times f(X)$ is a correspondence from $X$ to $f(X)$.
\end{exa}

The \textbf{composition} of two correspondences $C:X \rra Y$ and $D:Y \rra Z$ is the correspondence $ D \circ C:X \rra Z$ defined by
\[
D \circ C =\left\{ (x,z) \in X \times Z \ \middle| \ \exists y \in Y \ \text{such that} \ (x,y) \in C, \ (y,z) \in D \right\}.
\]
Let $C$ be a subset of $X \times Y$.
The \textbf{transpose} of $C$, denoted by $C^T$, is defined by
\[
C^T=\left\{(y,x) \in Y \times X \ \middle| \ (x,y) \in C\right\}.
\] 
A subset $C \subset X \times Y$ is a correspondence from $X$ to $Y$ if and only if its transpose $C^T$ is a correspondence from $Y$ to $X$.  
For a subset $C$ of $X \times Y$ and a subset $\sigma$ of $X$, its \textbf{image} $C(\sigma)$ is the set $\left\{y \in Y \ \middle| \ (x,y) \in C , \ x \in \sigma\right\}$.

When $(X,d_X)$ and $(Y,d_Y)$ are metric spaces, the \textbf{distortion} of a correspondence $C:X \rra Y$ is defined as follows:
\[
{\rm dis} \ C=\sup \left\{ | d_X(x,x')-d_Y(y,y') | \ \Big| \ (x,y),(x',y') \in C \right\}.
\]
The distortion of $C:X \rra Y$ is related to the distance between $X$ and $Y$ as follows.

\begin{prop}[Theorem 7.2.5 of \cite{BBI}]
\label{GHco}
For any two metric spaces $(X,d_X)$ and $(Y,d_Y)$,  
\[
d_{GH}(X,Y) =\frac{1}{2} \inf \left\{ {\rm dis} \ C \ \middle| \ C:X \rra Y \right\}, 
\]
where $d_{GH}$ is the Gromov-Hausdorff distance.
\end{prop}

We here recall the definitions of the Hausdorff distance and the Gromov-Hausdorff distance.
A \textbf{semi-metric} (psuedometric) is a function $d:X \times X \ra \R$ satisfying for all $x,y,z \in X$
\begin{enumerate}
\item $d(x,y) \geq 0$,
\item $d(x,x)=0$,
\item $d(x,y)=d(y,x)$,
\item $d(x,y)+d(y,z) \geq d(x,z)$.
\end{enumerate}
When $X$ and $Y$ are subspaces of some semi-metric space $(Z,d)$,
\[
d_H(X,Y):=\max \left\{\sup_{x \in X} \inf_{y \in Y} d(x,y), \sup_{y \in Y} \inf_{x \in X} d(x,y) \right\}
\]
is the \textbf{Hausdorff distance} between $X$ and $Y$.
When $(X,d_X)$ and $(Y,d_Y)$ are metric spaces with possibly different metrics,
it is possible to make a semi-metric $d$ on $X \sqcup Y$ satisfying $d|_{X \times X} =d_X$ and $d|_{Y \times Y}=d_Y$ (see \cite[Remark 7.3.12]{BBI} ).
Regarding $X$ and $Y$ as subspaces of $(X \sqcup Y,d)$, the Hausdorff distance between $X$ and $Y$ is thus defined.
We also denote the Hausdorff distance by $d_H(X,Y;d)$ to emphasize the dependency on the semi-metric.
Then the \textbf{Gromov-Hausdorff distance} is defined by 
\[
d_{GH}(X,Y):=\inf_d d_H(X,Y;d),
\]
where the infimum is taken over all semi-metrics on $X \sqcup Y$ satisfying $d|_{X \times X} =d_X$ and $d|_{Y \times Y}=d_Y$. 
Hence, Proposition \ref{GHco} provides a characterization of $d_{GH}(X,Y)$ by means of the distortions.

\subsection{Persistence modules}

A \textbf{persistence module} $\lV=\left( V_a,v^b_a \right)$ over the real numbers $\R$ is an indexed family of vector spaces
$\left\{V_a \ \middle| \ a \in \R\right\}$ and a doubly-indexed family of linear maps 
\[
\left\{v^b_a:V_a \ra V_b \ \middle| \ a \leq b \right\}
\]
satisfying for all $a \leq b \leq c$
\[
v^c_a=v^c_b \circ v^b_a , \ v^a_a=\one_{V_a}.
\]

Let $\lU=\left( U_a,u^b_a \right)$ and $\lV=\left( V_a,v^b_a \right)$ be persistence modules over $\R$.
Given $\ee > 0$, a \textbf{homomorphism of degree $\ee$} is a family of linear maps $\Phi=\left\{\phi_a :U_a \ra V_{a+\ee}\ \middle| \ a \in \R \right\}$ such that
$v^{b+\ee}_{a+\ee} \circ \phi_a=\phi_b \circ u^b_a$ whenever $a \leq b$.
We denote the set of homomorphisms of degree $\ee$ from $\lU$ to $\lV$ by ${\rm Hom}^{\ee}(\lU,\lV)$.
The \textbf{composition} of two homomorphisms $\Phi=\left\{\phi_a \right\} \in \Hom^{\ee}(\lU,\lV)$ and $\Psi=\left\{\psi_a\right\} \in \Hom^{\delta}(\lV,\lW)$ is the homomorphism of degree $\ee+\delta$ defined by:
\[
\Psi \circ \Phi =\left\{\psi_{a+\ee} \circ \phi_a:U_a \ra W_{a+\ee+\delta} \right\} \in \Hom^{\ee+\delta}(\lU,\lW).
\]
For any persistence module $\lV=\left( V_a,v^b_a \right)$,
\[
\one^{\ee}_{\lV}:=\left\{ v^{a+\ee}_a:V_a \ra V_{a+\ee} \ \middle| \ a \in \R \right\}
\]
is a homomorphism of degree $\ee$ from $\lV$ to $\lV$.

A \textbf{filtered simplicial complex} $\lS=\left\{S_a \ \middle| \ a \in \R \right\}$ is a family of simplicial complexes
such that $S_a$ is a subcomplex of $S_b$ whenever $a \leq b$.
When the vertex set is unchanged for all $S_a$ (say $X$), $X$ is called the vertex set of $\lS$.

\begin{exa}
Let $\lS=\left\{S_a \right\}$ be a filtered simplicial complex and $u^b_a:H_k(S_a) \ra H_k(S_b)$ be the linear map induced by the inclusion $S_a \hra S_b$ for all $a \leq b$.
Then $H_k(\lS):=\left( H_k(S_a),u^b_a \right)$ forms a persistence module.
\end{exa}

Let $\lS=\left\{S_a \right\},\lT=\left\{T_a \right\}$ be two filtered simplicial complexes with vertex sets $X,Y$, respectively.
A correspondence $C:X \rra Y$ is \textbf{$\ee$-simplicial} from $\lS$ to $\lT$ if, for any $a \in \R$ and any simplex $\sigma \in S_a$ (recall that $\sigma$ is a subset of $X$),
every finite subset of $C(\sigma)$ is a simplex of $T_{a+\ee}$.

Two persistence modules $\lU$ and $\lV$ are said to be \textbf{$\ee$-interleaved} if there exist two homomorphisms
\[
\Phi \in {\rm Hom}^{\ee}(\lU,\lV) , \ \Psi \in {\rm Hom}^{\ee}(\lV,\lU)
\]
such that $\Psi \circ \Phi = \one^{2\ee}_{\lU}$ and $\Phi \circ \Psi = \one^{2\ee}_{\lV}$;
in other words, if the following four diagrams commute whenever $a \leq b$:
$$
\xymatrix{
U_a \ar[dr]^{\phi_a} \ar[rr]^{u^b_a} & & U_b \ar[dr]^{\phi_b} &  & U_a \ar[dr]^{\phi_a} \ar[rr]^{u^{a+2 \ee}_a} & & U_{a+2\ee} \\
& V_{a+\ee} \ar[rr]^{v^{b+\ee}_{a+\ee}} & & V_{b+\ee} ,&  & V_{a+\ee} \ar[ur]^{\psi_{a+\ee}} & , 
}
$$

$$
\xymatrix{
& U_{a+\ee} \ar[rr]^{u^{b+\ee}_{a+\ee}} & & U_{b+\ee}  &  & U_{a+\ee} \ar[dr]^{\phi_{a+\ee}}  &  \\
V_a \ar[rr]^{v^b_a} \ar[ur]^{\psi_a} & & V_b \ar[ur]^{\psi_b} & , & V_a \ar[ur]^{\psi_a} \ar[rr]^{v^{a+2\ee}_a} & & V_{a+2\ee}  . 
}
$$

Interleaving is useful to study some features of $\lV$ by $\lU$ (and vice versa).
For instance, when $u^{a+2\ee}_a$ is nonzero, we can deduce that the vector space $V_{a+\ee}$ is nonzero (see the upper right diagram).

\subsection{Rips complex}
\begin{figure}[htbp]
 \begin{minipage}{0.48\hsize}
 \begin{center}
  \includegraphics[width=60mm]{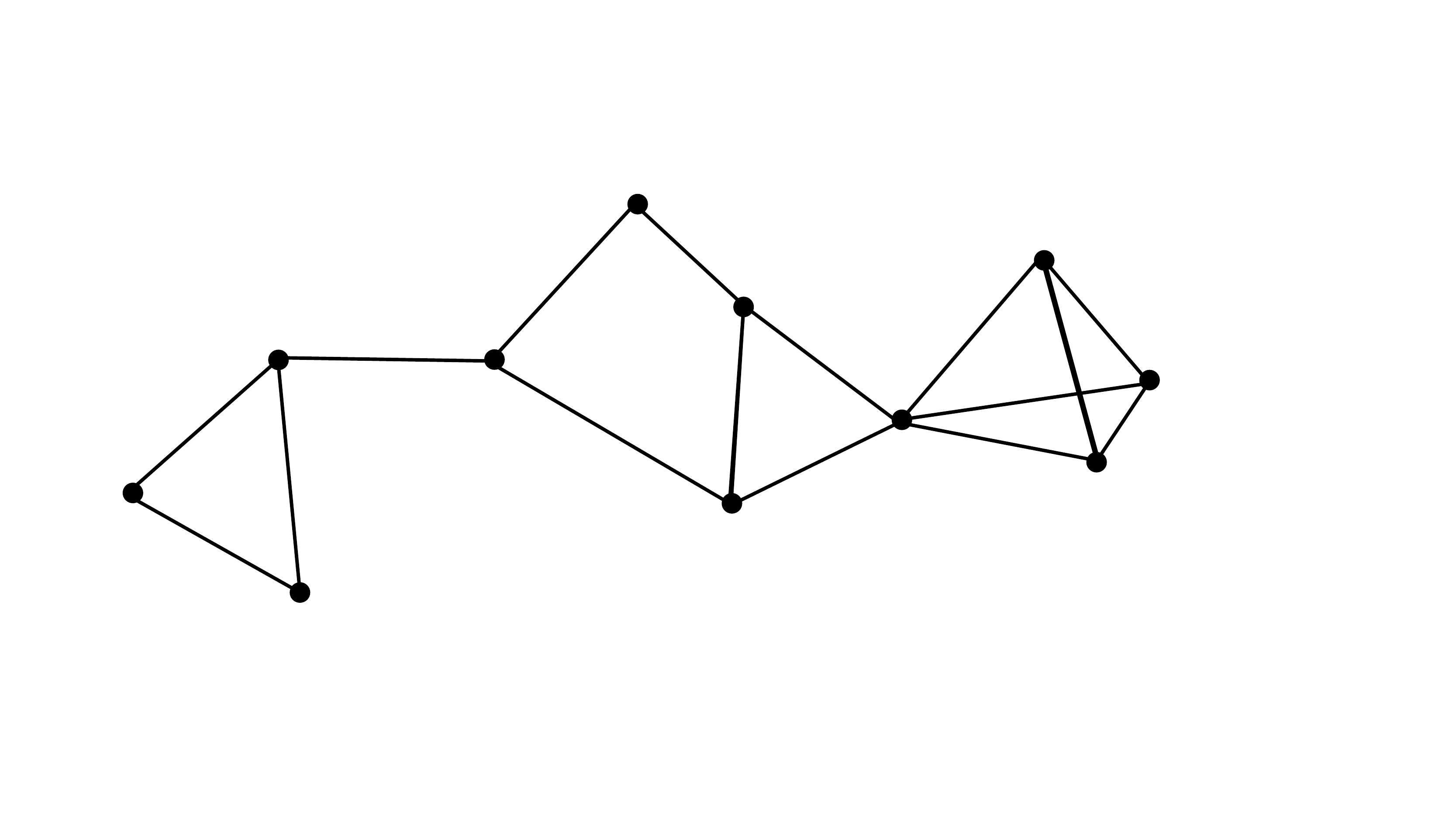}
 \end{center}
  \caption{A communication graph.}
  \label{fig:network}
 \end{minipage}
  \hspace{4mm}
 \begin{minipage}{0.48\hsize}
 \begin{center}
  \includegraphics[width=60mm]{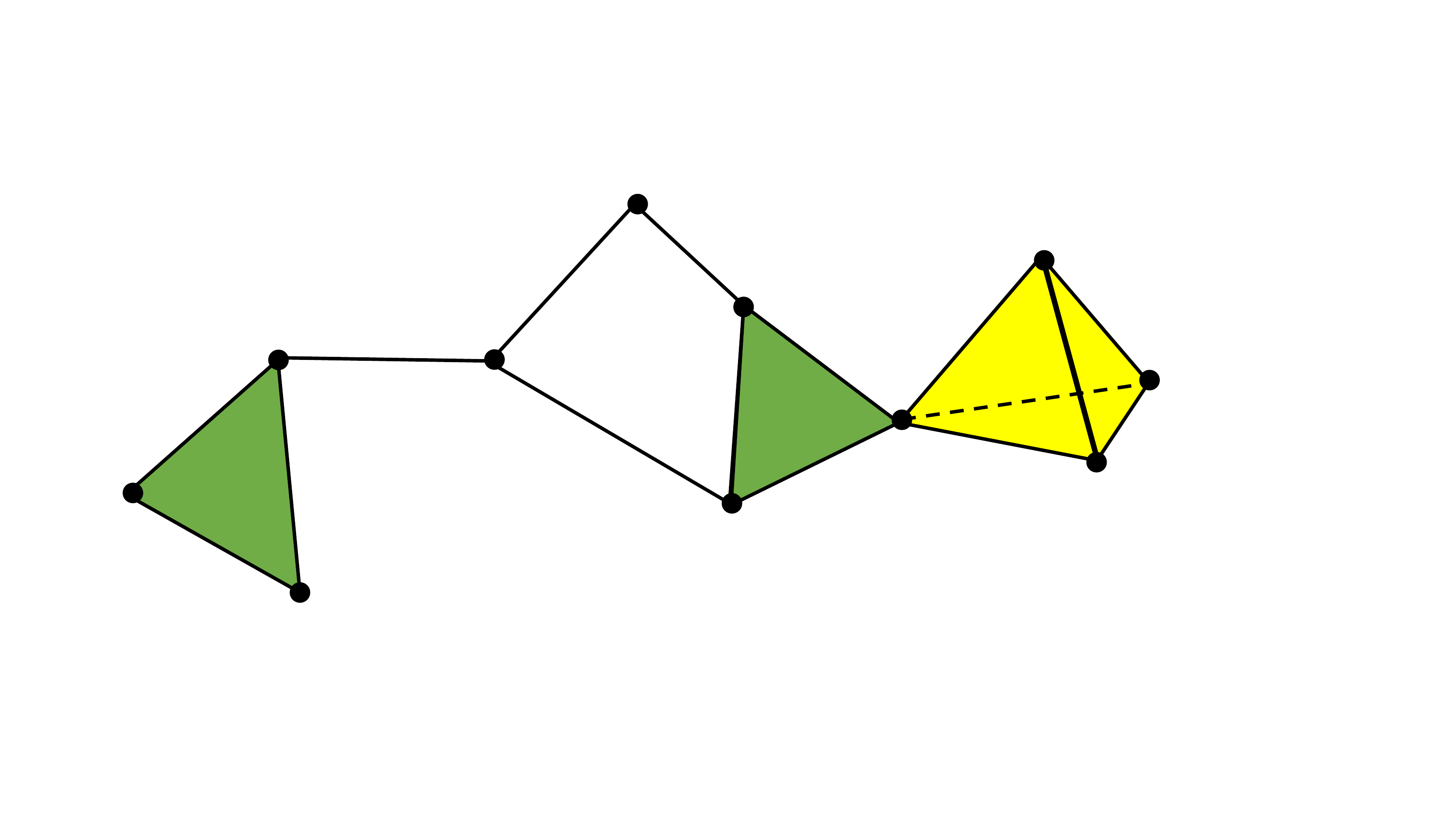}
 \end{center}
  \caption{The Rips complex of Figure \ref{fig:network}.}
  \label{fig:rips}
 \end{minipage}
\end{figure}

Let $(X,d_X)$ be a metric space.
For $a \geq 0$, we define a $k$-simplex $[x_{i_0} \cdots x_{i_k}]$ as a subset $\left\{ x_{i_0}, \ldots , x_{i_k} \right\}$ of $X$ that satisfies $x_{i_p} \neq x_{i_q} \ (p \neq q)$ and $d_X(x_{i_p},x_{i_q}) \leq a$ for all $p,q=0,\ldots ,k$. 
The set of these simplices forms a simplicial complex, called the \textbf{Rips complex} of $X$ with parameter $a$, denoted by $\cR (X;a)$.
For $a < 0$, we define $\cR(X;a)$ as a simplicial complex only consisting of the vertex set $X$ for convenience.
Since there is a natural inclusion $\cR(X;a) \hra \cR(X;b)$ whenever $a \leq b$,
$\cR(X)=\left\{\cR(X;a) \ \middle| \ a \in \R \right\}$ is a filtered simplicial complex on the vertex set $X$.
Figure \ref{fig:network} shows that the communication graph whose edges $[x_{i_p}x_{i_q}]$ are defined by $d(x_{i_p},x_{i_q}) \leq a$ and Figure \ref{fig:rips} shows its Rips complex.

\begin{prop}[Lemma 4.3 of \cite{CSO}]
\label{GHrips}
Let $(X,d_X)$ and $(Y,d_Y)$ be metric spaces.
If $d_{GH}(X,Y) < \frac{\ee}{2}$, then the persistence modules $H_k(\cR(X))$ and $H_k(\cR(Y))$ are $\ee$-interleaved.
\end{prop}

Proposition \ref{GHrips} is called the \textbf{Rips interleaving}.
The Rips interleaving is useful to study homological features of $\cR(Y)$ by $\cR(X)$, without directly computing the persistence module of $\cR(Y)$.

Using Proposition \ref{GHco}, we rewrite Proposition \ref{GHrips} as follows.

\begin{prop}
\label{corips}
Let $(X,d_X)$ and $(Y,d_Y)$ be metric spaces.
If there exists a correspondence $C:X \rra Y$ satisfying $\dis \ C \leq \ee$, then the persistence modules $H_k(\cR(X))$ and $H_k(\cR(Y))$ are $\ee$-interleaved.
\end{prop}
We extend Proposition \ref{corips} to a relative homology version in the next section.

\subsection{Optimal homologous cycles}

Here, we let the coefficient group be $\lF=\lR$.
Let $K$ be a simplicial complex of dimension greater than $p$.
Let $\{\sigma_j\}^{n}_{j=1}$ and $\{\tau_i\}^{m}_{i=1}$ be the sets of oriented $(p+1)$- and $p$-simplices in $K$, respectively.
We view a $p$-chain $x = \sum_i x_i \tau_i $ as a vector $x=(x_1,\cdots,x_m) \in \lR^m$.
In this way the $l^1$-norm of $x \in C_p(K)$ is defined by $\norm{x}_1=\sum_i \abs{x_i}$. 
Fixing bases of $C_{p+1}(K)$ and $C_p(K)$, we view the boundary map $\pa : C_{p+1}(K) \ra C_p(K)$ as a matrix $B$.

For a $p$-cycle $x$, the \textbf{optimal homologous cycle problem} \cite{DHK} is to solve:
\begin{equation}
\label{l1 chain problem}
\min \|\hat{x}\|_1 ~~ \mbox{ subject to  } \hat{x}= x+B y , ~ y \in \lR^n.
\end{equation}

When $L$ is a subcomplex of $K$ and $x + C_p(L) $ is a relative cycle in $C_p(K,L)$,  we choose a representative $x \in C_p(K)$ composed only by $p$-simplices in $K-L$.
Let $\{\sigma_j\}^{n}_{j=1}$ and $\{\tau_i\}^{m}_{i=1}$ be as above and the first $t$ and $s$ simplices of $\{\sigma_j\}^{n}_{j=1}$ and $\{\tau_i\}^{m}_{i=1}$ be all oriented $(p+1)$- and $p$-simplices in $L$, respectively.
Then the \textbf{optimal homologous relative cycle problem} is to solve:
\begin{align}
\label{l1 relative chain problem}
\min & ~ \|\hat{x}\|_1 \\
\mbox{subject to} &~  \hat{x}= x+B \ (0,\cdots,0,y_{t+1},\cdots,y_n)^T +(a_1,\ldots,a_s,0,\ldots,0) ^T, \nonumber \\ 
& y_j,a_i \in \lR. \nonumber
\end{align}
The vectors $(0,\cdots,0,y_{t+1},\cdots,y_n)$ and $(a_1,\ldots,a_s,0,\ldots,0)$ represent $(p+1)$ and $p$-chains in $K-L$ and $L$, respectively.

\section{RELATIVE INTERLEAVING}
\label{relative}
	Let $A$ be a subspace of a metric space $X$.
Then $\cR(A;a)$ is a subcomplex of $\cR(X;a)$ for all $a \in \R$.
We consider a pair $\cR(X,A)=(\cR(X),\cR(A))$ of filtered Rips complexes and its relative homology group $H_k(\cR(X,A;a))=H_k(\cR(X;a),\cR(A;a))$.
In this section, we extend some results in \cite{CSO} to relative versions, especially relative Rips interleaving.

\subsection{Relative correspondence}
	Let $X$ and $Y$ be two sets and let $A$ and $B$ be subsets of $X$ and $Y$, respectively.
A \textbf{correspondence} $C$ from $(X,A)$ to $(Y,B)$, denoted by $C:(X,A) \rra (Y,B)$, is a correspondence $C:X \rra Y$ satisfying $C(A) \subset B$ and $C^T(B) \subset A$.
If $C:(X,A) \rra (Y,B)$ is a correspondence, then the transpose $C^T$ is also a correspondence from $(Y,B)$ to $(X,A)$.

Given a map $f:X \ra Y$ satisfying $f(A) \subset B$, we write $f:(X,A) \ra (Y,B)$.
For $f:(X,A) \ra (Y,B)$, its graph $G(f)=\left\{(x,f(x)) \in X \times Y \ \middle| \ x \in X\right\} \subset X \times f(X)$ is not always a correspondence from $(X,A)$ to $(f(X),f(A))$ unlike Example \ref{graph}.
We can easily make counterexamples showing $G(f)^T(f(A)) \not\subset A$, because $G(f)^T(f(A))=f^{-1}(f(A))$. 
The following examples give some conditions on $f$ so that $f$ will satisfy $G(f)^T(f(A))=f^{-1}(f(A)) \subset A$.

\begin{exa}
If $f$ is injective, then $G(f)$ is a correspondence from $(X,A)$ to $(f(X),f(A))$. 
\end{exa}

\begin{exa}
If $f(A) \subset B$ and $f(X-A) \subset Y-B$, then $G(f)$ is a correspondence from $(X,A)$ to $(f(X),f(A))$.
\end{exa}

For a map $f:(X,A) \ra (Y,B)$, if $G(f)$ is a subset of $C:(X,A) \rra (Y,B)$, then the map $f:(X,A) \ra (Y,B)$ is said to be \textbf{subordinate} to $C$, denoted by $f:(X,A) \xrightarrow{C} (Y,B)$.

Let $A$ be a subset of $X$ and $\lS=\left\{S_a \right\}$ be a filtered simplicial complex with vertex set $X$.
The \textbf{restriction} $\lS^A=\left\{S^A_a \ \middle| \ a \in \R \right\}$ of $\lS$ to $A$ is the filtered simplicial complex such that $S^A_a$ is the maximal subcomplex of $S_a$ whose vertex set is $A$.
Each inclusion map $S_a \hra S_b$ induces a homomorphism $u^b_a:H_k(S_a,S^A_a) \ra H_k(S_b,S^A_b)$ whenever $a \leq b$,
and $H_k(\lS,\lS^A)=(H_k(S_a,S^A_a) , u^b_a)$ forms a persistence module.

Let $\lS=\left\{S_a\right\}$ and $\lT=\left\{T_a\right\}$ be two filtered simplicial complexes with vertex sets $X$ and $Y$, and let $A,B$ be subsets of $X,Y$, respectively.
A correspondence $C:(X,A) \rra (Y,B)$ is \textbf{$\ee$-simplicial} from $(\lS,\lS^A)$ to $(\lT,\lT^B)$ if $C:X \rra Y$ is $\ee$-simplicial, and for any $a \in \R$ and any simplex $\sigma \in S^A_a$,
every finite subset of $C(\sigma)$ is a simplex of $T^B_{a+\ee}$.

\begin{prop}
\label{unique}
Let $\lS=\left\{S_a\right\}$ and $\lT=\left\{T_a\right\}$ be two filtered simplicial complexes with vertex sets $X$ and $Y$, and let $A$ and $B$ be subsets of $X$ and $Y$, respectively.
Let $C:(X,A) \rra (Y,B)$ be $\ee$-simplicial from $(\lS,\lS^A)$ to $(\lT,\lT^B)$.
Then any subordinate map $f:(X,A) \xrightarrow{C} (Y,B)$ induces a canonical map $H_k(f) \in \Hom^{\ee}(H_k(\lS,\lS^A),H_k(\lT,\lT^B))$.
Moreover any maps induced by maps subordinate to $C$ are equal.
\end{prop}

\begin{proof}
Let $f:(X,A) \xrightarrow{C} (Y,B)$ be a subordinate map and $\sigma$ be a simplex in $S^A_a$.
Then $f(\sigma)$ is a finite subset of $C(\sigma)$.
This means $f(\sigma)$ is a simplex of $T^B_{a+\ee}$.
By applying the same argument to $X$ and $Y$, $f$ induces a simplicial map $(S_a,S^A_a) \ra (T_{a+\ee},T^B_{a+\ee})$ for all $a \in \R$.
Moreover the following diagram
$$
\xymatrix@C=0.2em{
(S_a,S^A_a)  \ar[dr] \ar[rr] & & (S_b,S^A_b) \ar[dr] & \\
& (T_{a+\ee},T^B_{a+\ee}) \ar[rr] & & (T_{b+\ee},T^B_{b+\ee}) 
}
$$
commutes whenever $a \leq b$, where the horizontal maps are inclusions.
Therefore $f$ induces $H_k(f) \in \Hom^{\ee} (H_k(\lS,\lS^A),H_k(\lT,\lT^B))$.

Any two subordinate maps $f_1,f_2:(X,A) \xrightarrow{C} (Y,B)$ induce simplicial maps from $(S_a,S^A_a)$ to $(T_{a+\ee},T^B_{a+\ee})$ that are contiguous (see \cite{Mu} for the definition and properties of contiguous maps).
In fact, for any $\sigma \in S^A_a$, the vertices of two simplices $f_1(\sigma)$ and $f_2(\sigma)$ span a simplex of $T^B_{a+\ee}$,
since these vertices comprise a finite subset of $C(\sigma)$.
Because any two contiguous maps are homotopic, we conclude that $H_k(f_1)=H_k(f_2)$. 
\end{proof}

Therefore the map $H_k(C):H_k(\lS,\lS^A) \ra H_k(\lT,\lT^B)$ is well-defined as $H_k(f)$ by a subordinate map $f:(X,A) \xrightarrow{C} (Y,B)$.

\begin{prop}
\label{inter}
Let $\lS$ and $\lT$ be two filtered simplicial complexes with vertex sets $X$ and $Y$, and let $A$ and $B$ be subsets of $X$ and $Y$, respectively.
If $C:(X,A) \rra (Y,B)$ is a correspondence such that $C$ and $C^T$ are both $\ee$-simplicial,
then $H_k(\lS,\lS^A)$ and $H_k(\lT,\lT^B)$ are $\ee$-interleaved.
\end{prop}

\begin{proof}
The diagonal set $\one_X:=\left\{(x,x) \ \middle| \ x \in X \right\}$ is a correspondence from $(X,A)$ to $(X,A)$ and satisfies $\one_X \subset C^T \circ C$.
Moreover, $\one_X$ is $2 \ee$-simplicial from $(\lS,\lS^A)$ to $(\lS,\lS^A)$.
The identity map from $(X,A)$ to $(X,A)$ is subordinate to $\one_X$ and it is also subordinate to $C^T \circ C$.
Thus by using Proposition \ref{unique}
\[
\one^{2 \ee}_{H_k(\lS,\lS^A)}=H_k(\one_X)=H_k(C^T \circ C) = H_k(C^T) \circ H_k(C).
\]
Similarly, $\one^{2 \ee}_{H_k(\lT,\lT^B)}=H_k(C) \circ H_k(C^T)$, and the proof is complete. 
\end{proof}

\subsection{Relative Rips interleaving}

If we find a correspondence between two filtered simplicial complexes satisfying the assumptions of Proposition \ref{inter},
then their persistence modules are interleaved.
Here, we construct such a correspondence for filtered Rips complexes.

\begin{thm}[Relative Rips interleaving]
\label{corelative}
Let $(X,d_X)$ and $(Y,d_Y)$ be metric spaces and let $A$ and $B$ be subspaces of $X$ and $Y$, respectively.
If there exists a correspondence $C:(X,A) \rra (Y,B)$ satisfying $\dis \ C \leq \ee$, then the persistence modules $H_k(\cR(X,A))$ and $H_k(\cR(Y,B))$ are $\ee$-interleaved.
\end{thm}

\begin{proof}
Let $\sigma$ be a finite subset of $A$ and $\tau$ be any finite subset of $C(\sigma)$.
For any $y,y' \in \tau$, there exist vertices $x,x'$ of $\sigma$ such that $(x,y),(x',y') \in C$.
If $\sigma$ is a simplex of $\cR(A;a)$, then $d_X(x,x') \leq a$ for any two vertices $x,x'$ of $\sigma$.
It follows from $\dis \ C \leq \ee$ that
\[
d_Y(y,y') \leq d_X(x,x') + \ee \leq a + \ee.
\]
Thus $\tau$ is a simplex of $\cR(B;a+\ee)$.
By substituting $X$ and $Y$ for $A$ and $B$ in the above argument, it can be seen that $C$ is $\ee$-simplicial from $\cR(X,A)$ to $\cR(Y,B)$. 
Symmetrically, $C^T$ is also $\ee$-simplicial from $\cR(Y,B)$ to $\cR(X,A)$.
Then the conclusion follows from Proposition \ref{inter}. 
\end{proof}

The following proposition is a simple extension of Proposition \ref{GHco}.

\begin{prop}
\label{reGHco}
If a correspondence $C:(X,A) \rra (Y,B)$ satisfies $\dis \ C \leq \ee$, then $d_{GH}(X,Y) \leq \frac{\ee}{2}$ and $d_{GH}(A,B) \leq \frac{\ee}{2}$.
\end{prop}

Conversely, with some additional assumptions, we show that $d_H(X,Y) \leq \frac{\ee}{2}$ and $d_H(A,B) \leq \frac{\ee}{2}$ imply the existence of $C:(X,A) \rra (Y,B)$ satisfying $\dis \ C \leq \ee$. 

\begin{prop}
\label{Hrelative}
Let $X_1,X_2,Y_1$, and $Y_2$ be subspaces of a metric space $(Z,d)$. 
If $d_H(X_i,Y_i) \leq \frac{\ee}{2}$ for $i=1,2$ and $X_1 \cap X_2 = Y_1 \cap Y_2 = \emptyset$,
then there exists a correspondence $C:(X_1 \sqcup X_2, X_j) \rra (Y_1 \sqcup Y_2, Y_j)$ satisfying $\dis \ C \leq \ee$ for $j=1,2$.
\end{prop}

\begin{proof}
If $d_H(X_i,Y_i) \leq \frac{\ee}{2}$, then the set
\[
C_i=\left\{ (x,y) \in X_i \times Y_i \ \middle| \ d(x,y) \leq \frac{\ee}{2} \right\}
\]
forms a correspondence $C_i:X_i \rra Y_i$.

Furthermore, we have $\dis \ C_i \leq \ee$.
Indeed, for any $x,x',y,y' \in Z$, the triangle inequalities
\begin{align*}
d(x,x') &\leq d(x,y)+d(y,x'),\ d(x',y) \leq d(x',y')+d(y,y') 
\end{align*}
imply
\[
d(x,x')-d(y,y') \leq d(x,y)+d(x',y').
\]
This means 
\[
\dis \ C_i \leq 2 \sup \{d(x,y) \mid (x,y) \in C_i \} \leq \ee.
\]

We write $X_{12}=X_1 \sqcup X_2$ and $Y_{12}=Y_1 \sqcup Y_2$.
Then the set
\[
C=\left\{(x,y) \in X_{12} \times Y_{12} \ \middle| \ (x,y) \in C_1 \ {\rm or} \ (x,y) \in C_2\right\}
\]
is a correspondence from $X_{12}$ to $Y_{12}$ satisfying $\dis \ C \leq \ee$. 

Moreover, from $X_1 \cap X_2 = Y_1 \cap Y_2 = \emptyset$, we have
\begin{align*}
C(X_j)
&= \left\{y \in Y_{12} \ \middle| \ (x,y) \in C, \ x \in X_j \right\} \\
&= \left\{y \in Y_{12} \ \middle| \ (x,y) \in C_j \right\} \\
&= \left\{y \in Y_j \ \middle| \ (x,y) \in C_j \right\} \\
& \subset Y_j.
\end{align*}
Similarly we obtain $C^T(Y_j) \subset X_j$, and hence $C$ is a correspondence from $(X_{12},X_j)$ to $(Y_{12},Y_j)$ satisfying $\dis \ C \leq \ee$. 
\end{proof}

As a corollary of Theorem \ref{corelative} and Proposition \ref{Hrelative}, we have the following.

\begin{cor}
\label{GHrelative}
Let $X_1,X_2,Y_1$, and $Y_2$ be subspaces of a metric space $(Z,d)$. 
If $d_H(X_i,Y_i) \leq \frac{\ee}{2}$ for $i=1,2$ and $X_1 \cap X_2 = Y_1 \cap Y_2 = \emptyset$,
then $H_k(\cR(X_1 \sqcup X_2,X_j))$ and $H_k(\cR(Y_1 \sqcup Y_2,Y_j))$ are $\ee$-interleaved for $j=1,2$.
\end{cor}

For $A \subset X \subset D \subset \R^d$, if $f:X \ra D$ satisfies $\sup_{x \in X} \|x-f(x)\| \leq \frac{\ee}{2}$, then $d_H(A,f(A)) \leq \frac{\ee}{2}$ and $d_H(X-A,f(X-A)) \leq \frac{\ee}{2}$.
Let us substitute $A,X-A,f(A)$, and $f(X-A)$ for $X_1,X_2,Y_1$, and $Y_2$, respectively.
If $f(A) \cap f(X-A)=\emptyset$, then we get the following.

\begin{cor}
\label{maprelative}
Let $A \subset X \subset D \subset \R^d$ and $f:X \ra D$ satisfy $\sup_{x \in X} \|x-f(x)\| \leq \frac{\ee}{2}$.
If $f(A) \cap f(X-A) = \emptyset$, then $H_k(\cR(X,A))$ and $H_k(\cR(f(X),f(A)))$ are $\ee$-interleaved.
\end{cor}

\section{REDUCING REDUNDANCIES}
\label{power}

The conclusion of the dSG criterion is that $D-\Nr \subset \cU(X)$.
However, it is possible that, for some $A \subset X$, the coverage $\cU(A)$ still contains $D-\Nr$.
From the point of view of power conservation of active sensors, we want to find the smallest subset $A$ whose coverage contains $D-\Nr$.
The topic about turing off active sensors is discussed in \cite{dSG1} for sensor networks in controlled boundaries.
In this section, we first review their works.
Then, we show the uniqueness of the generators achieving the coverage in the homological sense.
This property is applied into eliminating the redundant cover using the optimal homologous cycles.

A key concept in the proof of the dSG criterion is a realization map $\sigma: \cR(X;r_s) \ra \lR^d$, which sends a $p$-simplex $[x_0 \cdots x_p] \in \cR(X;r_s)$ to the singular $p$-simplex given by the convex hull of $x_0,\ldots,x_p$.
Under the assumptions $A1-A6$, $\sigma (\cR(F;r_s) )$ is in $\overline{\Nr}:=\lR^d-(D-\Nr)$ \cite[Lemma 3.1]{dSG2}.
Hence $\sigma$ induces the homomorphism 
\[
\sigma_*:H_d(\cR(X,F;r_s)) \longrightarrow H_d(\lR^d,\overline{\Nr}).
\]
From the naturality of long exact sequences, the following diagram
$$
\begin{CD}
H_d(\cR(X,F;r_s)) @>{\delta_*}>> H_{d-1}(\cR(F;r_s)) \\
@V_{\sigma_*}VV @V_{\sigma_*}VV \\
H_d(\lR^d,\overline{\Nr}) @>{\delta'_*}>> H_{d-1}(\overline{\Nr})
\end{CD}
$$
commutes, where $\delta_*$ and $\delta'_*$ are the connecting morphisms of the long exact sequences of $H_d(\cR(X,F;r_s))$ and $H_d(\lR^d,\overline{\Nr})$, respectively.

In fact, for any $\alpha \in H_d(\cR(X,F;r_s))$ satisfying $\iota_* \alpha \neq 0$,
it follows from the proof of dSG criterion that $\sigma_* \delta_* \alpha \neq 0$.
By commutativity, we get $\sigma_* \alpha \neq 0$.

Here, for a chain $z=\sum a_i \sigma_i$ ($a_i \neq 0$), the \textbf{coverage} of $z$ is defined to be 
\[
\cU(z)=\bigcup \{B_{r_c}(x) \mid x \mbox{ is a vertex of a simplex } \sigma_i, ~ a_i \neq 0\}.
\]

Then we get the following property by straightforward application of \cite[Corollary 4.1]{dSG1} to the dSG criterion.

\begin{prop}
For a relative cycle $z + C_d(\cR(F;r_s)) \in Z_d(\cR(X,F;r_s))$ such that $\sigma_* [z+ C_d(\cR(F;r_s))] \neq 0$, the coverage of $z$ contains $D-\Nr$.
\end{prop}


For this reason, we shall call $z \in C_d(\cR(X;r_s))$ satisfying $\sigma_* [z + C_d(\cR(F;r_s))] \neq 0$ a \textbf{coverage cycle} in $\cR(X,F)$.
In fact, all coverage cycles are homologous, by the following lemma.

\begin{lem}
\label{isomorphic}
For any coverage cycle $z$ in $\cR(X,F)$, $H_d(\cR(X,F;r_s))$ and $\lF [z] \oplus \ker \sigma_*$ are isomorphic.
\end{lem}

\begin{proof}
\label{unique generator}
The exact sequence
\[
0 \ra \ker \sigma_* \ra H_d(\cR(X,F;r_s)) \ra \im \sigma_* \ra 0
\]
is split because $\im \sigma_* $ is a vector space.
Thus $H_d(\cR(X,F;r_s))$ and $\im \sigma_* \oplus \ker \sigma_*$ are isomorphic.
The Alexander duality leads to 
\[
H_d(\lR^d,\overline{\Nr}) \simeq \lF
\]
because $D-\Nr$ is connected.
This implies that $\im \sigma_* \simeq \lF$ since $\sigma_*$ is a nonzero, which completes the proof.
\end{proof}

Let the coefficient group of homology groups be $\lF = \lR$.
For a coverage cycle $z$ in $\cR(X,F)$, let us set a solution $z^*$ of the problem \eqref{l1 relative chain problem}.
The solution $z^*$ is homologous to $z$ as a relative chain, and thus $z^*$ is also a coverage cycle in $\cR(X,F)$.

Let us call $\tilde{z}$ a \textbf{minimal} coverage cycle in $\cR(X,F)$ if $\tilde{z}$ is smaller than or equal to any coverage cycle in $\cR(X,F)$ with respect to $l^1$-norm.
Then the important point of Lemma \ref{isomorphic} is that any coverage cycle is homologous to a minimal coverage cycle in $\cR(X,F)$ as a relative chain.
Hence once we get a coverage cycle, a solution of the problem \eqref{l1 relative chain problem} is always a minimal coverage cycle in $\cR(X,F)$.

\section{COVERAGE OF PERTURBED SENSORS}
\label{perturbed sensors}

In this section, we consider the case where perturbation affects sensors.
A map $f:X \ra D$ is seen as perturbation, that is, for any sensor $x \in X$, $f(x) \in D$ is the sensor after the perturbation.

\subsection{Proof of the main theorem}
\begin{proof}
Since $\ f(F) \subset N_f(\pa D)$ and $ f(X-F) \subset D-N_f (\pa D)$, we have $f(F)=f(X) \cap N_f(\pa D)$.
Thus $f(F)$ is a set of fence sensors of $f(X)$.
Moreover 
\[
f(F) \cap f(X-F) \subset N_f(\pa D) \cap (D-N_f(\pa D)) = \emptyset,
\]
and it follows from Corollary \ref{maprelative} that $H_d(\cR(X,F))$ and $H_d(\cR(f(X),f(F)))$ are $\ee$-interleaved. 
Then the diagram

\begin{xy}
(0,15) *{H_d(\cR(X,F;r_s-\ee))}="A",
(50,15)*{ H_d(\cR(X,F;r_w-\ee))}="B",
(100,15)*{H_d(\cR(X,F;r_w+\ee))}="C",
(25,0)*{H_d(\cR(f(X),f(F);r_s))}="D",
(75,0)*{H_d(\cR(f(X),f(F);r_w))}="E",

\ar^{i_*} (18,15);(32,15)
\ar^{H_d(G(f))} (4,12);(20,3)
\ar^{H_d(G(f))} (54,12);(70,3)
\ar^{j_*} (43,0);(57,0)
\ar^{i'_*} (68,15);(82,15)
\ar^{H_d(G(f)^T)}  (80,3);(96,12)
\end{xy}

commutes, where all horizontal maps are linear maps induced by inclusions.
Because $\iota^{\ee}_*$ is nonzero and 
\begin{align*}
\iota^{\ee}_*
&= i'_* \circ i_* \\
&= H_d(G(f)^T) \circ H_d(G(f)) \circ i_* \\
&= H_d(G(f)^T) \circ j_* \circ H_d(G(f)),
\end{align*}
$j_*$ is nonzero.
Therefore it follows from the dSG criterion that $\cU(f(X))$ contains the restricted domain $D-\Nr$. 
\end{proof}

\subsection{Reducing redundancies after perturbation}

The main theorem yields information about not only the coverage but also generators by interleaving.
For instance, if $[z + C_d(\cR(F;r_s - \ee))] \in H_d(\cR(X,F,r_s-\ee))$ satisfies $\iota^{\ee}_*[z+C_d(\cR(F;r_s-\ee))] \neq 0$, we can automatically conclude that $f(z)$ satisfies $j_*[f(z)+C_d(f(F);r_s)] \neq 0$, which implies $f(z)$ is a coverage cycle in $\cR(f(X),f(F))$.

In section \ref{power}, we have discussed power conservation of active sensors.
From Lemma \ref{isomorphic}, the coverage cycle $f(z)$ is homologous to a minimal coverage cycle in $\cR(f(X),f(F))$ as a relative chain.
To find a minimal coverage cycle in $\cR(f(X),f(F))$, we consider a optimal homologous relative cycle problem \eqref{l1 relative chain problem}.

Now let $C_{d+1}(\cR(X;r_s))$ and $C_d(\cR(X;r_s))$ be viewed as $\lR^n$ and $\lR^m$, and the first $t$ and $s$ components of vectors in $\lR^n$ and $\lR^m$ be viewed in $(d+1)$- and $d$-simplices of $\cR(F;r_s)$, respectively.  
A solution of the following problem is a minimal coverage cycle in $\cR(X,F)$:
\begin{align}
\label{l1 relative chain problem perturb}
\min & ~ \| w \|_1 \\
\mbox{subject to} &~  w= f(z)+B \ (0,\cdots,0,y_{t+1},\cdots,y_n)^T +(a_1,\ldots,a_s,0,\ldots,0) ^T, \nonumber \\ 
& ~y_j,a_i \in \lR, \nonumber
\end{align}
where $B$ is a matrix representation of the boundary map 
\[
\pa:C_{d+1}(\cR(f(X);r_s)) \ra C_d(\cR(f(X);r_s)).
\]

Then, what we need is to construct $B$ and solve the problem \eqref{l1 relative chain problem perturb}.
We note that we do not need to compute $H_d(\cR(f(X),f(F);r_s))$ over again to obtain a minimal coverage cycle.

%
%

\section{CONCLUSION}

In this paper, we have presented a new coverage criterion stable under perturbation.
The main mathematical contribution of this paper is the extension of Rips interleaving to the relative version.
Here, relative Rips interleaving enables us to study the coverage of slightly moved sensors.
With our criterion, even if we only have the communication data of the sensors before perturbation,
we can determine the coverage of the perturbed sensors.
Moreover, by computing a generator for our criterion, we have discussed how to eliminate redundant cover for power conservation of perturbed sensors.

\subsection*{ACKNOWLEDGEMENT}

The authors wish to express their sincere gratitude to Emerson Escolar for valuable discussions on this paper.
The authors also thank an anonymous referee for the careful reading and variable comments.
This work is partially supported by JSPS 24684007.

%
%
%
%


\begin{thebibliography}{99}

\bibitem{BBI}
D. Burago, Y. Burago, and S. Ivanov. A Course in Metric Geometry. volume 33 of Graduate Studies in Mathematics. American Mathematical Society, Providence, RI. 2001.
\bibitem{CCGGY}
F. Chazal, D. Cohen-Steiner, M. Glisse, L.J. Guibas, and S.Y. Oudot. Proximity of persistence modules and their diagrams. In SCG. 237-246. 2009.
\bibitem{CSO}
F. Chazal, V. de Silva, and S. Oudot. Persistence stability for geometric complexes. Geometriae Dedicata. 2013.
\bibitem{DHK}
T.K. Dey, A.N. Hirani, B. Krishnamoorthy. Optimal homologous cycles, total unimodularity, and linear programming, SIAM J.Comput. 40(4) 1026-1044. 2011. 
\bibitem{dSG1}
V. de Silva and R. Ghrist. Coordinate-free coverage in sensor networks with controlled boundaries. Int. J. Robotics Research 25 1205-1222. 2006.
\bibitem{dSG2}
V. de Silva and R. Ghrist. Coverage in sensor networks via persistent homology. Algebraic and Geometric Topology 7 339-358. 2007.
\bibitem{EH}
H. Edelsbrunner and J. Harer. Computational Topology: an Introduction. American Mathematical Society, Providence, RI. 2010.
\bibitem{EHi}
E.G. Escolar, Y. Hiraoka. "Computing Optimal Cycles of Homology Groups", A Mathematical Approach to Research Problems of Science and Technology. R. Nishii et al. ed. Springer Mathematics for Industry 5 101-118. 2014.
\bibitem{Ha}
A. Hatcher. Algebraic Topology, Cambridge Univ Press. 2001.
\bibitem{MA}
R. Mulligan and H. M. Ammari. Coverage in Wireless Sensor Networks: A Survey. Network Protocols and Algorithms vol. 2, 27-53. 2010.
\bibitem{MK}
S. Meguerdichian, F. Koushanfar, M. Potkonjak, and M. Srivastava. Coverage problems in wireless ad-hoc sensor network. IEEE INFOCOM vol.3, 1380-1387. 2001
\bibitem{Mu}
J. R. Munkres. Elements of Algebraic Topology. Westview Press. 1984.
\end{thebibliography}
\end{document}